\documentclass[a4paper,english]{lipics-v2018}

\usepackage{stackengine}
\usepackage{tikz-cd}
\usepackage[inference]{semantic}
\usepackage{proof}
\usepackage{float}
\usepackage{varwidth}
\usepackage{bussproofs}
\DeclareMathOperator{\type}{\text{ type}}
\DeclareMathOperator{\linear}{\text{ linear}}
\DeclareMathOperator{\diag}{\mathbf{Diag}}

\tikzset{%
    symbol/.style={%
        draw=none,
        every to/.append style={%
            edge node={node [sloped, allow upside down, auto=false]{$#1$}}}
    }
}
 \newtheorem{proposition}[theorem]{Proposition}

\usepackage{microtype}
\title{A diagram model of linear dependent type theory}

\titlerunning{Models of linear dependent type theory}

\author{Martin Lundfall}{Stockholm University}{martin@dapp.org}{}{}

\authorrunning{M. Lundfall}

\Copyright{Martin Lundfall}
\subjclass{
;
}

\category{}

\relatedversion{\\
  \url{http://kurser.math.su.se/pluginfile.php/16103/mod_folder/content/0/2017/2017_47_report.pdf}}

\supplement{}

\funding{}

\acknowledgements{I want to thank Peter Lumsdaine for his guidance and inspiration, without which this work would not have been possible.}

\EventEditors{}
\EventNoEds{0}
\EventLongTitle{Preprint}
\EventShortTitle{Preprint}
\EventAcronym{}
\EventYear{2018}
\EventDate{}
\EventLocation{}
\EventLogo{}
\SeriesVolume{}
\ArticleNo{}
\nolinenumbers 
\hideLIPIcs  

\begin{document}

\maketitle

\begin{abstract}
We present a type theory dealing with non-linear, ``ordinary'' dependent types (which we will call \textit{cartesian}), and \textit{linear types}, where both constructs may depend on terms of the former. In the interplay between these, we find the new type formers $\sqcap_{x : A}B$ and $\sqsubset_{x : A}B$, akin to $\Pi$ and $\Sigma$, but where the dependent type $B$, (and therefore the resulting construct) is a linear type. These can be seen as internalizing universal and existential quantification of linear predicates. We also consider two modalities, $M$ and $L$, transforming linear types into cartesian types and vice versa.
The theory is interpreted in a split comprehension category \cite{jacobs} $\pi : \mathcal{T} \to \mathcal{C}^\to$ accompanied by a split symmetric monoidal fibration, $\pi: \mathcal{L} \to \mathcal{C}$. This structure determines, for any context $\Gamma \in \mathcal{C}$, fibers $\mathcal{T}_\Gamma$ and $\mathcal{L}_\Gamma$, which become the category of cartesian types and the monoidal category of linear types over $\Gamma$, respectively. Here, the type formers $\sqcap_{x:A}$ and $\sqsubset_{x:A}$ are understood as right and left adjoints of the monoidal reindexing functor $\pi_A^* : \mathcal{L}_\Gamma \to \mathcal{L}_{\Gamma.A}$. The operators $M$ and $L$ give rise to a fiberwise adjunction $L \dashv M$ between $\mathcal{L}$ and $\mathcal{T}$, where the traditional exponential modality is understood as the comonad $! = LM$.

We provide a model of this theory called the \textit{Diagram model}, which extends the groupoid model of dependent type theory \cite{hofmann1998} to accommodate linear types. Here, cartesian types are interpreted as a family of groupoids, while linear types are interpreted as diagrams $A:\Gamma\to\mathcal{V}$ in any symmetric monoidal category $\mathcal{V}$. We show that the diagrams model can under certain conditions support a linear analogue of the univalence axiom, and provide some discussion on the higher-dimensional nature of linear dependent types.
\end{abstract}
\newpage
\tableofcontents
\section{Introduction \& summary of results}
Lately, there has been an increasing interest in combining linear and dependent types  \cite{schreiber2014quantization}, \cite{krishnaswami}, \cite{vakar14}, \cite{nothing}. The idea is that such a theory would inherit the higher-order nature of dependent types, while maintaining a careful account of how assumptions are used in a derivation. It is not completely clear, however, what the synthesis looks like, since in dependent type theory, variables may appear in both terms and types, but linear type theory only allows each variable to appear freely exactly once. Here, we take an approach inspired by \cite{krishnaswami} and \cite{vakar14}, in which we distinguish between non-linear, dependent types (which we call \textit{cartesian}), and linear types, and circumvent the issue by only allowing cartesian terms to appear in types (both cartesian and linear).

The theory splits contexts into two parts, divided by a semicolon, where the first part contains cartesian assumptions, for which weakening and contraction is admissible, while the second part contains linear assumptions, for which only exchange is allowed. We introduce two new type formers, $\sqcap_{x : A}B$ and $\sqsubset_{x : A}B$, akin to $\Pi$ and $\Sigma$, but where the dependent type $B$ (and therefore the resulting construct) is a linear. The traditional $!$ modality is deconstructed as a comonad arising from the adjoint pair $L \dashv M$, where $L$ is a functor (or modality) sending cartesian types into linear, and $M$ sends linear types to cartesian. We have $\Pi_{x : A}B_M \cong (\sqcap_{x :A}B)_M$, for linear $B$, and, assuming a few additional rules, a linear isomorphism $(\Sigma_{x :A}C)_L \cong \sqsubset_{x :A}C_L$ for cartesian $C$.

Compared to ordinary dependent type theory, we get additional elimination and computation rules for both $\Sigma$ and $Id$-types when eliminating into a linear type.

We postulate the existence of two universes, $L$ and $U$, containing codes of linear and cartesian types, respectively and assumed to be closed under all type formers.

We develop categorical semantics for the theory by defining a model as a \textit{comprehension category} \cite{jacobs}, $\pi : \mathcal{T} \to \mathcal{C}$ equipped with a \textit{split symmetric monoidal fibration} $q : \mathcal{L} \to \mathcal{C}$ over the same base. A split symmetric monoidal fibration has just enough structure to make the fibers $\mathcal{L}_\Gamma$ over a context $\Gamma \in \mathcal{C}$ into symmetric monoidal categories, and reindexing functors (strict) monoidal functors. The traditional linear type formers $\&, \oplus, 0, \top, \multimap$ correspond to the existence of binary products and coproducts, initial and terminal object and internal homs in each fiber, such that these are preserved under reindexing. The new type formers $\sqcap_{x : A}B$ and $\sqsubset_{x :A}B$ correspond to right and left adjoints to the reindexing functor $\pi_A^* : \mathcal{L}_{\Gamma} \to \mathcal{L}_{\Gamma.A}$, while the modalities $L$ and $M$ give rise to a fiber adjunction between $\mathcal{L}$ and $\mathcal{T}$. The new rules for $\Sigma$ are automatically satisfied by the semantic interpretation of $\Sigma_A$ as a left adjoint to the reindexing functor $\pi_A^* : \mathcal{T}_\Gamma \to \mathcal{T}_{\Gamma.A}$. The new rules for $Id$-types impose an additional condition on the semantic interpretation of $Id$, which are always fulfilled if our identity types are extensional.

We consider two concrete models, the first being the \textit{families model}, in which cartesian types consist of families of sets, indexed by their context set $\Gamma$, and a linear type in the context $\Gamma$ is a $\Gamma$-indexed family of objects in a given symmetric monoidal category $\mathcal{V}$. Examples of suitable $\mathcal{V}$ supporting all type formers present in our syntax are $\mathbf{AbGrp}$, $\mathbf{GCTop}_*$, $\mathbf{Vect}_F$, i.e. the category of abelian groups, the category of compact generated, pointed topological spaces and the category of vector spaces over a field $F$, respectively.

Generalizing the families model, we get the \textit{diagrams model}, in which contexts are interpreted as groupoids, and cartesian types over a groupoid $\Gamma$ are diagrams in $\mathbf{Gpd}$ over $\Gamma$, and linear types over $\Gamma$ are diagrams in a given symmetric monoidal category $\mathcal{V}$ over $\Gamma$. Just as the groupoid model \cite{hofmann1998} can be shown to support a univalent universe, we construct a linear analogue of the univalence axiom and show that it holds in the diagrams model if the adjunction $L \dashv M$ factors through sets.
\section{Syntax}\label{syntax}
As cartesian type formers, we use the standard $\Sigma$, $\Pi$, and identity type formers as well as universe types $U$ and $L$ for linear and cartesian types, respectively.
The purely linear part of our type theory contains all the type formers of intuitionistic linear logic; the additive connectives $\&,\oplus,\mathbf{0},\top$ and the multiplicatives $\otimes,\mathbf{1},\multimap$. In addition to these, we have the new type formers $\sqsubset$, $\sqcap$, which play a role analogous to that of $\Sigma$ and $\Pi$ in the cartesian setting.
Finally, we have the two modalities, $M$ and $L$, which turns linear types into cartesian, and vice versa.
A detailed presentation of our syntax can be found in \cite{MLundfall}. For the familiar, ``purely'' dependent or linear type formers, our presentation offers no significant surprises, except for a couple of additional rules for $\Sigma$ and the identity type. Therefore, we focus on presenting the syntax for the new type formers $\sqsubset$, $\sqcap$ and the modalities $M$ and $L$.

\subsection{Auxiliary elimination rules}\label{newElim}
Besides the traditional rules for $\Pi$, $\Sigma$ and the identity type, we find that since we can now eliminate into linear types, we must introduce an extra elimination and computational rule for each one. These additional rules are presented next along with the traditional elimination rules in Figure \ref{elimfig}.
\begin{figure}
\fbox{
  \begin{minipage}{.40\linewidth}
      \footnotesize
  \[
\inference[{\footnotesize$\Sigma$-E$_1$}]{\Gamma, t : \Sigma_{x : A}B \vdash C \type\\
\Gamma, x : A, y : B \vdash c : C[(x, y)/t]\\
\Gamma \vdash s : \Sigma_{x : A}B}
{\Gamma \vdash \hat c[s] : C[s/t]}
\]\\
\[
\inference[{\footnotesize =-E$_1$}]{\Gamma, x, y : A, p : x =_Ay \vdash C \type\\
\Gamma, z : A \vdash c: C[z/x, z/y, \text{refl}(z)/p]\\
\Gamma \vdash M, N : A\\
\Gamma \vdash P : M =_A N}
{\Gamma \vdash R^{Id} : C[M/x, N/y, P/p]}
\]
\end{minipage}
  \begin{minipage}{.60\linewidth}
    \footnotesize
\[
  \inference[{\footnotesize $\Sigma$-E$_2$}]{\Gamma, t : \Sigma_{x : A}B \vdash C \linear\\
\Gamma, x : A, y : B; \Xi \vdash c : C[(x, y)/t]\\
\Gamma \vdash s : \Sigma_{x : A}B}
{\Gamma; \Xi[pr_1(s)/x][pr_2(s)/y] \vdash \hat c[s] : C[s/t]}
\]\\
\[
  \inference[{\footnotesize =-E$_2$}]{\vdash \Gamma, x, y : A, p : x =_Ay; \Xi \text{ ctxt}\\
    \Gamma, x, y : A, p : x =_Ay \vdash C \linear\\
\Gamma, z : A; \Xi[z/x, z/y, \text{refl}(z)/p] \vdash c: C[z/x, z/y, \text{refl}(z)/p]\\
\Gamma \vdash M, N : A\\
\Gamma \vdash P : M =_A N}
{\Gamma; \Xi[M/x, N/y, P/p] \vdash R^{Id} : C[M/x, N/y, P/p]}
\]
\end{minipage}
}
\caption{Elimination rules for $\Sigma$ and $Id$}
\label{elimfig}
\end{figure}

\subsection{The modalities $M$ and $L$}\label{syntLandM}
We introduce two the modal operators $M$ and $L$, which transfers a linear type/term to its cartesian counterpart and vice versa. Semantically, this will establish a fiberwise monoidal adjunction between the categories of linear and cartesian types:
\[
\begin{tikzcd}
\mathcal{L}_{\Gamma} \ar[r, "M"{name=A, below}, bend right] & \mathcal{T}_{\Gamma} \ar[l, "L"{name=B, above}, bend right] \ar[from=A, to=B, symbol=\vdash]
\end{tikzcd}
\]
where the exponential modality from traditional linear logic is understood as the comonad $! = LM$. The decomposition of the exponential into an adjunction goes back to at least \cite{benton1995mixed}, and is given an interesting new light in \cite{licata2017fibrational}, where it is seen as a particular case of a more general procedure of encoding structure in contexts.

The rules for the operators $M$ and $L$ are presented in Figure \ref{MLfig}.
\begin{figure}
\fbox{
  \small
\begin{minipage}{.6\linewidth}
\[
\inference{\Gamma \vdash A \type}{\Gamma \vdash A_L \linear}[L-F]
\]
\[
\inference{\Gamma \vdash a : A}{\Gamma ; \cdot \vdash a_L : A_L}[L-I]
\]
\[
\inference{(\Gamma \vdash B \linear) \\
(\vdash \Gamma; \Xi' \text{ ctxt})\\ \Gamma; \Xi \vdash y : A_L \quad \Gamma, x : A; \Xi' \vdash t : B}{\Gamma; \Xi, \Xi' \vdash \text{let  $x$ be $y$ in }t : B}[L-E]
\]
\[
\inference{\Gamma; \Xi \vdash \text{let  $x$ be $s_L$ in }t : B}{\Gamma; \Xi \vdash \text{let  $x$ be $s_L$ in }t \equiv t[s/x]  : B}[L-C]
\]
\[
  \inference{\Gamma; y : A_L, \Xi \vdash t : B \\
    \Gamma; \Xi' \vdash a : A_L 
  }
  {\Gamma; \Xi, \Xi' \vdash \text{let $x$ be $a$ in $t[x_L/y]$} \equiv t[a/y] : B}[L-U]
\]

\end{minipage}
\begin{minipage}{.4\linewidth}
  \small
\[
\inference{\Gamma \vdash B \linear}{\Gamma \vdash B_M \type}[M-F]
\]
\[
\inference{\Gamma ; \cdot \vdash b : B}{\Gamma \vdash \sigma(b) : B_M}[M-I]
\]
\[
\inference{
\Gamma \vdash t : B_M\\
}{\Gamma; \cdot \vdash \sigma^{-1}(t) : B \\}[M-E]
\]
\[
\inference{
\Gamma \vdash \sigma(b)  : B_M\\
}{\Gamma; \cdot \vdash \sigma^{-1}(\sigma(b)) \equiv b : B \\}[M-C$_1$]
\]
\[
\inference{
\Gamma; \cdot \vdash \sigma^{-1}(t) : B_M\\
}{\Gamma \vdash \sigma(\sigma^{-1}(t)) \equiv t : B_M \\}[M-C$_2$]
\]
\end{minipage}
}
\caption{Typing rules for $M$ and $L$}
\label{MLfig}
\end{figure}
The interpretation of $L$ and $M$ as an adjoint pair is already present at the syntactic level. We can show that they form instances of a Haskell-like \textbf{Functor} class, by constructing terms: $\textup{\texttt{fmapM}} : (A \multimap B)_M \to A_M \to B_M$ and $\textup{\texttt{fmapL}} : L(A \to B) \multimap (LA \multimap LB)$, satisfying the functor laws.

Furthermore, we can construct a ``counit'' $\epsilon : LM \implies 1$ satisfying the universal property of adjunction (thanks to L-U). The syntactic formulation of the statement becomes:
\begin{theorem}[$L \dashv M$]
There is a term $\Gamma; \beta_1 : B_{LM} \vdash \epsilon_B: B$ with the following property:\\
For any term: $\Gamma; y : A_L \vdash f : B$, there is a unique term $\Gamma, x : A \vdash g : B_{M}$ such that $\Gamma; y : A_L \vdash \epsilon_B[\text{let $x$ be $y$ in }g_L/\beta_1] \equiv f : B$.
\end{theorem}
 Based on this knowledge we expect the right adjoint $M$ to preserve limits, and indeed we find an isomorphism: $A_M \times B_M \cong (A \& B)_M$.
We can now also reformulate some common results about the exponential modality using $! = LM$, such as $(A \& B)_{LM} \cong A_{LM} \otimes B_{LM}$.
\subsection{$\sqcap$ and $\sqsubset$}
Since we allow linear types to depend on terms of cartesian types, we can form new versions of the $\Pi$- and $\Sigma$-types, denoted $\sqcap$ and $\sqsubset$, respectively. The typing rules for these are presented in Figure \ref{clawhousefig}.
\begin{figure}
\fbox{
  \begin{minipage}{.47\linewidth}
    \small
\[
\inference{\Gamma \vdash A \type \quad \Gamma, x : A \vdash B \linear}{\Gamma \vdash \sqcap_{x : A}B \linear}[$\sqcap$-F]
\]\\
\[
  \inference{\vdash \Gamma; \Xi \text{ ctxt}\\
    \Gamma, x : A; \Xi \vdash b : B
  }
{\Gamma; \Xi \vdash \lambda x. b : \sqcap_{x: A}B}[$\sqcap$-I]
\]\\
\[
\inference{\Gamma; \Xi \vdash t : \sqcap_{x : A}B \quad \Gamma \vdash a : A}
{\Gamma; \Xi \vdash t(a) : B[a/x]}[$\sqcap$-E]
\]\\
\[
\inference{\Gamma; \Xi \vdash \lambda x. b (a) : \sqcap_{x : A}B}
{\Gamma; \Xi \vdash \lambda x. b (a) \equiv b[a/x] : B[a/x]}[$\sqcap$-C]
\]
\end{minipage}
  \begin{minipage}{.53\linewidth}
    \small
\[
\inference{\Gamma \vdash A \type \quad \Gamma, x : A \vdash B \linear}{\Gamma \vdash \sqsubset_{x : A}B \linear}[$\sqsubset$-F]
\]\\
\[
\inference{\Gamma \vdash s : A \quad \Gamma; \Xi \vdash b : B[s/x]}{\Gamma; \Xi \vdash (s, b) : \sqsubset_{x : A}B}[$\sqsubset$-I]
\]\\
\[
  \inference{\vdash \Gamma;\Xi' \text{ ctxt}\\
    \Gamma, x : A \vdash C \linear\\
    \Gamma; \Xi \vdash t : \sqsubset_{x : A}B \quad \Gamma, x : A; \Xi', y : B \vdash c : C}{\Gamma; \Xi, \Xi' \vdash \text{let $x, y$ be $t$ in c : C}}[$\sqsubset$-E]
\]\\
\[
\inference{\Gamma; \Xi \vdash \text{let $x, y$ be $(s, t)$ in c : C}}
{\Gamma; \Xi \vdash \text{let $x, y$ be $(s, t)$ in c} \equiv c[s/x][t/y] : C}[$\sqsubset$-C]
\]
\end{minipage}
}
\caption{Typing rules for $\sqsubset$ and $\sqcap$}
\label{clawhousefig}
\end{figure}
The sense in which $\sqcap$ and $\sqsubset$ are ``linear analogues'' of $\Pi$ and $\Sigma$ can be formalized in the following way:
\begin{proposition}\label{M-sqcap}
  For all $\Gamma \vdash A \type$ and $\Gamma, x : A \vdash B \linear$, there is an isomorphism:
  \[
    \Pi_{x : A}B_M \cong (\sqcap_{x : A}B)_M
    \]
\end{proposition}

We would like to show a similar result relating $\Sigma$ and $\sqsubset$, but for this we need a couple of additional rules. First, we assume uniqueness rules for $\Sigma$ and $\sqsubset$, asserting that the elimination rule followed by the introduction rule is the identity. In other words, for any $p : \Sigma_{x : A}B$ and $q : \sqsubset_{x : C}D$, we have $(pr_1(p),pr_2(p)) \equiv p$ and $\text{let $x, y$ be $q$ in $(x, y)$} \equiv q$ \footnote{The former is provable as a propositional identity \cite[Corollary~2.7.3]{hott-book}. Perhaps it is possible to obtain a similar result for $\sqsubset$, using the ``surrogate equality'' described in the end of Section \ref{syntLandM}}:
Second, we assume a kind of naturality rule for the $L$ modality:
\[
  \inference{
    \Gamma; \Xi, y : B \vdash e : C\\
    \Gamma, x : A; \Xi'\vdash u : B\\
    \Gamma; \Xi'' \vdash t : A_L
}
{ \Gamma; \Xi, \Xi', \Xi'' \vdash e[\text{let $x$ be $t$ in $u$}/y] \equiv \text{let $x$ be $t$ in $e[u/y]$} : C}[Nat$_L$]
\]
\begin{proposition}\label{L-subset}
  Assuming the Nat$_L$ and the uniqueness rules for $\Sigma$ and $\sqsubset$, there is a linear isomorphism\footnote{Here a linear isomorphism, $A \cong B$, means a pair $f : A \multimap B$, $g : B \multimap A$ such that the composite is judgmentally equal to the identity. We discuss the weaker notion of linear equivalence in Section \ref{highermodel}}
  \[
    (\Sigma_{x :A}B)_L \cong \sqsubset_{x:A}B_L
  \]
\end{proposition}
As outlined in section \ref{semtypform}, the semantic interpretation of the type formers $\Pi$, $\sqcap$ and $\Sigma$, $\sqsubset$ are as right and left adjoints to reindexing functors respectively. Based on this interpretation we can understand these equivalence results through the diagram:
\[
\begin{tikzcd}[row sep=huge, column sep=large]
  \mathcal{L}_{\Gamma.A} \ar[d] \ar[r, "M_{\Gamma.A}"{name=C, above}, bend right]& \mathcal{T}_{\Gamma.A} \ar[l, "L_{\Gamma.A}"{name=D, below}, bend right, swap] \ar[from=D, to=C, symbol=\dashv] \ar[d]\\
  \mathcal{L}_\Gamma \ar[u, "\sqcap_A"{name=A}, bend right,swap] \ar[u, "\sqsubset_A"{name=B}, bend left] \ar[r, "M_{\Gamma}"{name=C, above}, bend right] \ar[from=B, to=A, symbol=\dashv,xshift=-5] \ar[from=B, to=A, symbol=\dashv,xshift=4]
  & \mathcal{T}_\Gamma \ar[l, "L_{\Gamma}"{name=D, below}, bend right,swap] \ar[u, "\Pi_A"{name=F}, bend right,swap] \ar[u, "\Sigma_A"{name=E}, bend left] \ar[from=D, to=C, symbol=\dashv] \ar[from=E, to=F, symbol=\dashv,xshift=-5] \ar[from=E, to=F, symbol=\dashv,xshift=4]
\end{tikzcd}
\]
  
\section{Semantics}\label{semantics}
\subsection{Structural semantic core}
Our semantic exploration of linear dependent type theory begins with the notion of a model. For the cartesian fragment of our theory, we follow \cite{jacobs} and ask for a \textit{comprehension category}, $\pi : \mathcal{T} \to \mathcal{C}^\to$, where $\mathcal{C}$ is a category of context with terminal object, and the fibrations $\mathcal{T}_\Gamma$ contains the cartesian types over $\Gamma$. For the linear fragment of our theory, we would like a fibration $q : \mathcal{L} \to \mathcal{C}$ where each fiber $\mathcal{L}_\Gamma$ is a symmetric monoidal category and the reindexing functors are symmetric monoidal. This is captured in the notion of a \textit{(lax) monoidal fibration}:
\begin{definition}
A \textbf{lax monoidal fibration} \cite{zawadowski} is a fibration $p : E \to B$ along with
\begin{enumerate}
\item Two functors $\otimes : E \times_{B} E \to E$ and $I : B \to E$ fitting into the following diagram:
\[
\begin{tikzcd}
E \times_{B} E \ar[r, "\otimes"] \ar[rd] & E \ar[d, "p"] & B \ar[l, "I",swap] \ar[ld, "1_B"] \\
& B &
\end{tikzcd}
\]
\item Three fibred natural isomorphisms $\alpha, \lambda$ and $\rho$ associated with the diagrams:
\[
\begin{tikzcd}
E \times_B E \times_B E \ar[r, "1_E \times_B \otimes"] \ar[d, "\otimes \times_B 1_E",swap] & E \times_B E \ar[d, "\otimes", swap] \\
E \times_B E \ar[r, "\otimes"] \ar[ru, "\alpha", Rightarrow, shorten <=10pt,shorten >=10pt]  & E
\end{tikzcd}
\]
and
\[
\begin{tikzcd}
B \times_B E \ar[r, "I \times_B 1_{E}"] \ar[rdd, "\pi_2",swap] & E \times_B E  \ar[dd, "\otimes",swap] & E \times_B B \ar[ldd, "\pi_1"] \ar[l, "1_E \times I",swap] \\
\ar[r, "\lambda", xshift=10 pt, Rightarrow, shorten <=20pt, shorten >=20pt] & {}& \ar[l, "\rho", Rightarrow, xshift=-10pt, shorten <=20pt, shorten >=20pt,swap] \\
  & E &  
\end{tikzcd}
\]
\item such that $\alpha$, $\lambda$ and $\rho$ satisfies the pentagon and triangle identities in each fiber.
\item for every $b \in B$, $\rho_{I_b} = \lambda^{-1}_{I_b} : I_b \otimes I_b \to I_b$
\end{enumerate}
\end{definition}
To avoid any coherence problems, we require the both the comprehension category and the monoidal fibration to be \textit{split}.
\begin{definition}
A \textbf{model for linear dependent type theory} consists of a split comprehension category $\pi : \mathcal{T} \to \mathcal{C}^\to$ and a split symmetric monoidal fibration $q : \mathcal{L} \to \mathcal{C}$, as illustrated in the following picture:
\[
\begin{tikzcd}
\mathcal{L} \ar[rd, "q"] & \mathcal{T} \ar[d, "p"] \ar[r, "\pi"] & \mathcal{C^\rightarrow} \ar[ld, "\text{cod}"]\\
& \mathcal{C}
\end{tikzcd}
\]
\end{definition}
where $\text{cod}$ denotes the codomain fibration functor.

This provides the necessary machinery to interpret all the structural rules of our theory as well as the rules for $\otimes$ and $I$, by constructing an interpretation function $[[-]]$, which sends:
  \begin{itemize}
  \item Cartesian contexts $\Gamma$ to objects of $\mathcal{C}$, considered up to judgmental equality and renaming of bound variables.
  \item Linear contexts $\Xi = a_1 : A_1, a_2 : A_2, \dots a_n : A_N$ in $\Gamma$ to objects $[[\Xi]] = \bigotimes^n_{i = 1}[[A_i]]$ of $\mathcal{L}_{[[\Gamma]]}$.
  \item Cartesian types $A$ in $\Gamma$ to objects of $\mathcal{T}_{[[\Gamma]]}$.
  \item Linear types $B$ in $\Gamma$ to objects of $\mathcal{L}_{[[\Gamma]]}$.
  \item Cartesian terms $M : A$ in $\Gamma$ to sections of the projection morphism $\pi([[A]]) : [[\Gamma,A]] \to [[\Gamma]]$.
  \item Linear terms $b : B$ in $\Gamma; \Xi$ to morphisms $[[b]] : [[\Xi]] \to [[B]]$.
  \end{itemize}
\subsection{Semantic type formers}\label{semtypform}
Equipped with the baseline structure of a model in which we can interpret the structural rules of our theory, we formulate the conditions under which such models support various type formers.
From now on, we will assume that the comprehension category comprising the core of our syntax is full, i.e. that the functor $\pi : \mathcal{T} \to \mathcal{C}^{\to}$ is full and faithful. This simplifies the semantic interpretation of many type formers.

The interpretation of the purely linear type formers $\otimes, I, \multimap, \&, \oplus, \top$ and $0$ in symmetric monoidal categories is well known. See for instance \cite{mellies}. Notice that $\otimes$ and $I$ types are supported in any model. For a model to support the type formers $\multimap, \&, \oplus, \top$ and $0$, correspond to the condition that the fibers of $\mathcal{L}$ have weak versions of internal homs, binary products and coproducts, and terminal and initial object, and that these are stable under reindexing functors.
\subsubsection{$\Pi$ and $\Sigma$}
What it means for a model of linear dependent type theory to support $\Pi$-types is directly inherited from the standard, non-linear case; we require right adjoints to reindexing functors satisfying a Beck-Chevalley condition.

As the rules $\Sigma$ contains one more eliminator than usual ($\Sigma$-E$_2$ from Section \ref{newElim}), one might wonder whether this poses additional clauses in the definition of the semantic type former. But as it turns out, the relevant condition will always hold in any model supporting $\Sigma$-types:
\begin{definition}A model of LDTT \textbf{supports $\Sigma$-types} if it satisfies the following:
  \begin{enumerate}
  \item For all $A \in \mathcal{T}_{\Gamma}$, the induced functor $\pi_A^* : \mathcal{T}_{\Gamma} \to \mathcal{T}_{\Gamma.A}$ has a left adjoint, $\Sigma_A$,
  \item (Beck-Chevalley) such that for all pullbacks
    \[
    \begin{tikzcd}
     \Gamma.E \ar[d, "\pi_E"] \ar[r, "q_{E, E'}"] & \Delta.E' \ar[d, "\pi_{E'}"] \\
     \Gamma \ar[r, "f"] & \Delta
    \end{tikzcd}
    \] 
    the natural transformation: $\Sigma_Eq^* \to f^*\Sigma_{E'}$ is a natural isomorphism, and
  \item the induced map $pair: \Gamma.A.B \to \Gamma.\Sigma_AB$ is an isomorphism
  \end{enumerate}
\end{definition}
We will denote the inverse of $pair$ by $(pr_1, pr_2)$, when it exists.
This structure is sufficient to support new elimination rule ($\Sigma$-E$_2$).
\begin{theorem}
If a model of LDTT supports $\Sigma$-types, then for every object $C \in \mathcal{L}_{\Gamma.\Sigma_AB}$, morphism $c : \Xi \to (pair_{A,B})^*C$ in $\mathcal{L}_{\Gamma.A.B}$ and section $s : \Gamma \to \Gamma.\Sigma_AB$, there exists a morphism $\hat c_s : s^*(pr_1, pr_2)^*\Xi \to s^*C$ such that given sections $a : \Gamma \to \Gamma.A$ and $b : \Gamma.A \to \Gamma.A.B$ we have $\hat c_{(a, b)} = a^*b^*c : a^*b^*\Xi \to a^*b^*C$.
\begin{proof}
Let $\hat c_s = c(pr_1,pr_2)^*s^*c$. Given sections $a : \Gamma \to \Gamma.A$ and $b : \Gamma.A \to \Gamma.A.B$ we compose with $pair$ to get the section $(a, b) = pair \circ ba : \Gamma \to \Gamma.\Sigma_AB$. We have:
    \[
     \hat c_{(a,b)} = (pr_1,pr_2)^*(a,b)^*c = (pr_1,pr_2)^*pair^*b^*a^*c = b^*a^*c
    \]
  \end{proof}
\end{theorem}
\subsubsection{Identity types}
The situation for Id-types requires a bit more care. If one wants to keep the theory intensional, we need to add condition (2) to make sure that the semantic identity types satisfy the added elimination rule, =-E$_2$.

\begin{definition}[Id-types]\label{idsemantic}
  A model of LDTT \textbf{supports Id-types} if, for all $A \in \mathcal{T}_{\Gamma}$, there exists an object $Id_A \in \mathcal{T}_{\Gamma.A.\pi_A^*A}$ and a morphism $r_A : \Gamma.A \to \Gamma.A.\pi_A^*A.Id_{A}$ such that $\pi_{Id_A} \circ r_A = v_A$, and:
  \begin{enumerate}
\item For any commutative diagram:
  \[
    \begin{tikzcd}
    \Gamma.A \ar[d, "r_A"] \ar[r] & \Delta.C \ar[d, "\pi_C"]\\
    \Gamma.A.\pi_A^*A.Id_A \ar[r] & \Delta
    \end{tikzcd}
  \]
  there exists a lift $J: \Gamma.A.\pi_A^*A.Id_A \to \Delta.C$ making the two triangles commute.
\item For any pair of objects, $C, \Xi \in \mathcal{L}_{\Gamma.A.\pi_A^*A.Id}$, sections $M, N : \Gamma \to \Gamma.A$, $P : \Gamma \to \Gamma.M^*(N^+)^*Id_A$, and morphism $c : r_A^*\Xi \to r_A^*C$, there exists a morphism $\hat c_{[M,N,P]} : M^*(N^+)^*P^*\Xi \to M^*(N^+)^*P^*C$ such that $\hat c_{[M,M,\text{refl}]} = M^*c$.
\end{enumerate}
\end{definition}
Notice that if our type theory has extensional id-types, in the sense that $a =_A b$ implies $a \equiv b$, then the second condition is always met.
\subsubsection{$\sqcap$- and $\sqsubset$-types}
The semantic type formers for the linear dependent $\sqcap$ and $\sqsubset$ is akin to that of $\Pi$ and $\Sigma$. They are given by adjoints to the functors between fibers of $\mathcal{L}$ induced by the projection maps in $\mathcal{C}$.
\begin{definition}
  A model of LDTT \textbf{supports $\sqcap$-types} if, for all $A \in \mathcal{T}_{\Gamma}$, the induced monoidal functor $\pi_A^* : \mathcal{L}_{\Gamma} \to \mathcal{L}_{\Gamma.A}$ has a monoidal right adjoint, $\sqcap_A$ satisfying the following Beck-Chevalley condition:

  For all pullback squares in $\mathcal{C}$ of the following form:
  \[
      \begin{tikzcd}
    \Gamma.E \ar[d, "\pi_E"] \ar[r, "q_{E, E'}"] & \Delta.E' \ar[d, "\pi_{E'}"] \\
    \Gamma \ar[r, "f"] & \Delta
    \end{tikzcd}
  \]
the canonical natural transformation $f^*\sqcap_{E'} \to \sqcap_{E}q^*_{E, E'}$ is a natural isomorphism.\\
\end{definition}
\begin{definition}
  It \textbf{supports $\sqsubset$-types} if, for all $A \in \mathcal{T}_{\Gamma}$, the functor every $\pi_A^*$ has a monoidal left adjoint, satisfying the following:
  \begin{enumerate}
  \item (Beck-Chevalley): For all pullbacks squares as above, the natural transformation $\sqsubset_Eq^* \to f^*\sqsubset_{E'}$ is a natural isomorphism.
  \item (Frobenius reciprocity): For all objects $\Xi \in \mathcal{L}_{\Gamma}$ and $B \in \mathcal{L}_{\Gamma.A}$, the canonical morphism $\sqsubset_A(\Xi\{\pi_A\} \otimes B) \to \Xi \otimes \sqsubset_AB$ is an isomorphism.
  \end{enumerate}
\end{definition}
\subsubsection{The operators $M$ and $L$}
\begin{definition}\label{semanticML}
  A model of LDTT with unit \textbf{supports the operators $M$ and $L$} if there exists functors $M : \mathcal{L} \leftrightarrow \mathcal{T} : L$ which are cartesian with respect to the fibrations $p : \mathcal{T} \to \mathcal{C}$ and $q : \mathcal{L} \to \mathcal{C}$, such that
  \begin{itemize}
  \item $L \dashv M$ is a fibred adjunction,
  \item $L(1) \cong I$
  \item and there is an isomorphism of hom-sets:
  \[
    \mathcal{L}_{\Gamma.A}(\pi_A^*(\Xi'), \pi_A^*(B)) \cong \mathcal{L}_\Gamma(LA \otimes \Xi', B)\].

  \end{itemize}
\end{definition}
Recall that a fibred adjunction implies that there are natural isomorphisms making the following diagram commute:
\[
  \begin{tikzcd}
  \mathcal{L}_{\Gamma.A}  \ar[r, "M_{\Gamma.A}", bend right] & \mathcal{T}_{\Gamma.A} \ar[l, "L_{\Gamma.A}", bend right] \\
  \mathcal{L}_{\Gamma} \ar[r, "M_{\Gamma}", bend right] \ar[u, "\pi^*_A"] & \mathcal{T}_{\Gamma} \ar[l, "L_{\Gamma}", bend right]  \ar[u, "\pi^*_A"] \\
  \end{tikzcd}
\]
which from a syntactic perspective ensures that $M$ and $L$ commute with substitution.

Note that the interpretation of a term $\Gamma \cdot \vdash \sigma(a) : A_M$ arises from the adjunction via $\sigma : \mathcal{L}_\Gamma(I, A) \cong \mathcal{L}_\Gamma(L(1), A) \cong \mathcal{T}_\Gamma(1, M(A))$.

The final condition of the definition is what yields the elimination and computation rules L-U, and while it might appear somewhat unnatural semantically, it does turn out to hold in a broad variety of models, due to the following result:
\begin{theorem}In a model of LDTT that supports $\multimap$ type formers, then any fibred adjunction $L \dashv M$ where $L(1) \cong I$ satisfies $\mathcal{L}_{\Gamma.A}(\pi^*_A(\Xi'), \pi^*(B)) \cong \mathcal{L}_{\Gamma}(LA \otimes \Xi', B)$.
  \begin{proof}
    A model supporting internal homs must have reindexing functions which preserve these. That is, we have an isomorphism $\pi_A^*[\Xi, B] \cong [\pi_A^*\Xi, \pi^*_AB]$. We get a chain of isomorphisms:
\[
\begin{split}
&\mathcal{L}_{\Gamma}(LA \otimes \Xi, B) \cong 
\mathcal{L}_{\Gamma}(LA, [\Xi, B]) \cong 
\mathcal{T}_{\Gamma}(A, M_{\Gamma}[\Xi, B]) \cong \\
&\mathcal{T}_{\Gamma.A}(1, \pi_A^*(M_{\Gamma}[\Xi, B])) \cong 
\mathcal{T}_{\Gamma.A}(1, M_{\Gamma.A}\pi_A^*[\Xi, B])) \cong 
\mathcal{L}_{\Gamma.A}(L_{\Gamma.A}(1), \pi_A^*[\Xi, B])) \cong \\
&\mathcal{L}_{\Gamma.A}(I, \pi_A^*[\Xi, B])) \cong 
\mathcal{L}_{\Gamma.A}(I, [\pi_A^*\Xi, \pi_A^*B])) \cong 
\mathcal{L}_{\Gamma.A}(\pi_A^*\Xi, \pi_A^*B).
\end{split}
\]
  \end{proof}
\end{theorem}
\section{Diagram Model}\label{models}
The main novelty of this paper is the Diagram model of linear dependent type theory. This model extends the groupoid model of dependent type theory \cite{hofmann1998} to support linear types, while still maintaining a higher dimensional interpretation of the identity type. Most interestingly, perhaps, it provides a model in which we can support univalent universes, both for cartesian and linear types. The diagram model can be seen as a natural generalization of the set indexed families model described by \cite{vakar14}. We briefly recall the set indexed families model below as a useful comparison to the diagrams model.
\subsection{Set indexed families}
\begin{definition}[$Fam(\mathcal{C})$]
For an arbitrary category $\mathcal{C}$, let $Fam(\mathcal{C})$ denote the category whose objects consists of pairs $(S, f)$ where $S$ is a set and $f$ is a function $f : S \to Ob(\mathcal{C})$. Morphisms of $Fam(\mathcal{C})$ are pairs $(u, \alpha) : (S, f) \to (S', g)$ where $u : S \to S'$ and $\alpha : S \to \text{Mor}(\mathcal{C})$ such that $\alpha(s) : f(s) \to g(u(s))$ for all $s \in S$.
\end{definition}
By projecting a family to its indexing set, we get a fibration $p : Fam(\mathcal{C}) \to \mathbf{Set}$ and a comprehension category by defining $\pi(S, f) = fst: \{(s, t) \; | \; s \in S, t : \top \to f(s)\} \to S$. \footnote{As long as $\mathcal{C}$ has a terminal object and the hom-sets $\mathcal{C}(\top,A)$ are small for all $A \in \mathcal{C}$.}

Letting $\mathcal{C} = \mathbf{Set}$ thus gives us a (full, split) comprehension category, forming the cartesian part of our model. For the linear part, we can for any symmetric monoidal category $\mathcal{V}$ form a monoidal fibration by a simple pointwise construction, giving us the following picture:
\[
\begin{tikzcd}
Fam(\mathcal{V}) \ar[rd, "q"]  & Fam(\mathbf{Set}) \ar[d, "p"] \ar[r, "\pi"] & \mathbf{Set}^\to \ar[ld, "cod"] \\
& \mathbf{Set}
\end{tikzcd}
\]

In this setting, most type formers will be given by a simple pointwise construction which are preserved under reindexing. It turns out that the families model supports the type formers $\otimes, I, \multimap, \oplus, 0, \&$, and $\top$ if $\mathcal{V}$ is a monoidal category which is closed, has binary coproducts, initial object, binary products and terminal object respectively.

It supports $\sqcap$-types if $\mathcal{V}$ has small products, and $\sqsubset$ if $\mathcal{V}$ has small coproducts that distribute over $\otimes$ (Frobenius reciprocity).

The families model of course also supports $\Pi$ and $\Sigma$-types, and since its identity types are extensional, the extra condition posed on our semantic identity types poses no additional difficulty.

Whenever $\mathcal{V}$ is a concrete category, the adjunction $F \dashv U$ will induce a fiber adjunction between the corresponding fibrations, which forms support for the operators $M$ and $L$ (as long as $F(\mathbf{1}) \cong I$).
\subsection{Diagrams in monoidal categories}
For any category $\mathcal{C}$, there is a fibration $cod : \diag(\mathcal{C}) \to \mathbf{Cat}$, where $\diag(\mathcal{C})$ refers to the category of diagrams in $\mathcal{C}$, i.e. consisting of objects $J : \mathcal{A} \to \mathcal{C}$, and morphisms between $J : \mathcal{A} \to \mathcal{C}$ and $J' :\mathcal{B} \to \mathcal{C}$ are functors $F : \mathcal{A} \to \mathcal{B}$ equipped with a natural transformation $\eta: J \implies J' \circ F$. In other words, the fibers of $\diag(\mathcal{C})$ are functor categories, which we write $[\Gamma, \mathcal{C}]$, for any small category $\Gamma$. Any functor $F : \mathcal{A} \to \mathcal{B}$ in the base induces a canonical lift $F^* : [\mathcal{B}, \mathcal{C}] \to [\mathcal{A},\mathcal{C}]$ simply given by precomposition.

When $\mathcal{C}$ has a terminal object $\top$ such that the collections $\mathcal{C}(\top, A)$ are small for any $A \in \mathcal{C}$, we form a comprehension category:
\[
\begin{tikzcd}
\diag(\mathcal{C}) \ar[d, "dom"] \ar[r,"\pi"] & \mathbf{Cat}^\to \ar[ld, "cod"] \\
\mathbf{Cat}
\end{tikzcd}
\]
where the functor $\pi$ sends a diagram $A : \Gamma \to \mathcal{C}$ to the \textbf{Grothendieck construction} for $A$, i.e. the category whose objects are pairs $(\gamma, t_\gamma)$ where $\gamma \in \Gamma$, $t_\gamma : \top \to A(\gamma)$. Morphisms $(\gamma, t_\gamma) \to (\gamma', t'_{\gamma'})$ consists of morphisms $u : \gamma \to \gamma'$ such that $ A(u) \circ t_\gamma = t'_{\gamma'}$. \footnote{If $\mathcal{C}$ is a 2-category, this can be weakened so that morphisms $(\gamma, t_\gamma) \to (\gamma', t'_{\gamma'})$ are pairs $(u, \alpha)$, where $u : \gamma \to \gamma'$ and $\alpha$ is a 2-cell $\alpha : A(u) \circ t_\gamma \implies t'_{\gamma'}$.}

When $\mathcal{C}$ is any symmetric monoidal category $\mathcal{V}$, there is an obvious symmetric monoidal structure on each fiber $[\Gamma, \mathcal{V}]$, given pointwise.

Restricting the base of the fibration to groupoids instead of categories, and setting $\mathcal{C} = \mathbf{Gpd}$ we get a model of linear dependent type theory which expands the groupoid model by Hofmann and Streicher \cite{hofmann1998}:
\[
\begin{tikzcd}
\diag(\mathcal{V}) \ar[rd, "dom"]  & \diag{\mathbf{(Gpd)}} \ar[d, "dom"] \ar[r, "\pi"] & \mathbf{Gpd}^\to \ar[ld, "cod"] \\
& \mathbf{Gpd}
\end{tikzcd}
\]
Since $\top$ is the groupoid $\mathbf{1}$ consisting of a single object, we will equate the functor $t_\gamma : 1 \to A(\gamma)$ with an object $a_\gamma$ of $A(\gamma)$, and the natural transformation $\alpha : A(u) \circ t_\gamma \implies t'_{\gamma'}$ with a morphism $\alpha_\gamma : A(u)(a_\gamma) \to a'_{\gamma'}$.

As shown in \cite{hofmann1998}, this model supports $\Pi$ and $\Sigma$ type formers, and provides an interesting interpretation of the identity type $Id_A$ as the arrow category of $A$.

This construction satisfies the additional requirement in our definition of semantic identity types:
\begin{theorem}Given a $\Gamma$-indexed groupoid $A$, diagrams $C, \Xi \in [\Gamma.A^\to, \mathcal{V}] \cong [\Gamma.A.\pi_A^*A,Id_A, \mathcal{V}]$, sections $M, N : \Gamma \to \Gamma.A$, $P : \Gamma.A \to \Gamma.Id_A$ and a natural transformation $c : \Xi \circ r_A \implies C\circ r_A$, there exists a natural transformation $\hat c_{[M,N,P]} : \Xi \circ P^+ \circ N^+ \circ M \implies C \circ P^+ \circ N^+ \circ M$ such that $\hat c_{[M,M,\text{refl}]} = c \circ M$ \footnote{Where the sections $N^+ : \Gamma.A \to \Gamma.A.\pi_A^*A$ and $P^+ : \Gamma.A.\pi_A^* \to \Gamma.A.\pi_A^*.Id_A$ are weakenings of $N$ and $P$, i.e. functors ignoring the additional arguments}
  \begin{proof}
    The key point to observe is that there is always an isomorphism $(\gamma, P_\gamma : M_\gamma \to N_\gamma) \cong (\gamma, 1_{M_\gamma} : M_\gamma \to M_\gamma)$ given by the commutative diagram:
    \[
      \begin{tikzcd}
      (\gamma, M_\gamma) \ar[d, "P"] \ar[r, "1_M"] & (\gamma, M_\gamma) \ar[d, "1_M"]\\
      (\gamma, N_\gamma) \ar[r, "P^{-1}"] & (\gamma, M_\gamma)
      \end{tikzcd}
    \]
    forming a collection of isomorphisms in $\Gamma.A^\to$ giving rise to a natural isomorphism $\phi : r_A \circ M \implies P^+ \circ N^+ \circ M$. We define $\hat c[M,N,P]$ as the composite:
\[
\Xi \circ P^+ \circ N^+ \circ M \xrightarrow{\Xi_{\phi}} \Xi \circ r_A \circ M \xrightarrow{cM} C \circ r_A \circ M \xrightarrow{C_{\phi^{-1}}} C \circ P^+ \circ N^+ \circ M
\]
To see that the computation rule holds, we only need to notice that when $M \equiv N$ and $P = \text{refl}(M)$, $\phi$ is the identity natural transformation.
\end{proof}
\end{theorem}

As in the families model, limits and colimits are constructed pointwise, and preserved by precomposition, so the model supports $\&$, $\top$, $\oplus$, $0$, if $\mathcal{V}$ has binary products, terminal object, binary coproducts and initial object respectively.

When it comes to $\multimap$, we utilize the following result:
\begin{theorem}\label{diaghoms}If $\mathcal{V}$ has internal homs and is complete, $[\mathcal{C},\mathcal{V}]$ also has internal homs, defined for $F, G \in [\mathcal{C}, \mathcal{V}]$ by the end:
    \[
    [F, G] := \int_{x \in \mathcal{C}}[Fx, Gx])
    \]
\end{theorem}
These are preserved under reindexing, implying that the diagrams model supports $\multimap$ if $\mathcal{V}$ is monoidal closed and complete.
\begin{definition}
  For any functor $p : \mathcal{A} \to \mathcal{B}$ in the base, a left or right adjoint to the induced functor $p^* : [\mathcal{B}, \mathcal{V}] \to [\mathcal{A}, \mathcal{V}]$ is called a \textbf{left or right Kan extension} along $p$.
\end{definition}
We recall the following fact about Kan extensions:
  \begin{theorem}\label{kanlimits}
    Left (right) Kan extensions along $p : \mathcal{A} \to \mathcal{B}$ between two arbitrary small categories $\mathcal{A}$ and $\mathcal{B}$ exists if and only if $\mathcal{V}$ has all colimits (limits).
  \end{theorem}
The result above ensures the existence of left and right adjoints to reindexing functors in the diagrams model as long as $\mathcal{V}$ is co-complete or complete, respectively. Since our reindexing functors are given by precomposition they will always satisfy the Beck-Chevalley condition.
Again, in order to support $M$ and $L$, we can lift an adjunction between $\mathcal{V}$ and $\mathbf{Gpd}$ to fiber adjunction between the respective diagram categories. Therefore, for any diagrams model which supports $\multimap$ to support $M$ and $L$, it suffices to display an adjunction
\[
  \begin{tikzcd}
\mathcal{V} \ar[r, "L_0"{name=D, below}, bend right]  & \ar[l, "M_0"{name=C, above}, bend right] \mathbf{Gpd}
  \end{tikzcd}
\]
such that $L(1) \cong I$.
\begin{remark}\label{factorsets}
When the functor $\mathcal{V}(I, -) : \mathcal{V} \to \mathbf{Set}$ has a left adjoint $F$, we get an adjunction between $\diag(\mathcal{V})$ and $\diag(\mathbf{Gpd})$, induced by:
\[
  \begin{tikzcd}
\mathcal{V} \ar[r, "{\mathcal{V}(I, -)}"{name=D, below}, bend right]  & \ar[l, "F"{name=C, above}, bend right] \ar[from=C, to=D, symbol=\dashv] \mathbf{Set} \ar[r, "\delta"{name=B, below}, bend right]   &  \mathbf{Gpd} \ar[l, "\pi_0"{name=A, above}, bend right] \ar[from=A, to=B, symbol=\dashv]
  \end{tikzcd}
\]
where $\pi_0$ is the functor sending a groupoid to its set of connected components.
\end{remark}
\begin{theorem}\label{M-faith}There are models in which $M$ is not faithful.
  \begin{proof}
    Let $\mathcal{V}$ to be $\mathbf{Gpd}$ so that $L = \delta \pi_0$ and $M = \delta \mathbf{Gpd}(1, -)$. This induces a fiber adjunction $L \dashv M$ where $L(1) = 1$, but $M$ is not faithful.
  \end{proof}
\end{theorem}
\subsubsection{Universes in the diagrams model}
To support universes, assuming one inaccessible cardinal allows us to shift our perspective to from the category of small groupoids, $\mathbf{Gpd}$, to the category $\mathbf{GPD}$ of all groupoids. Among the objects of $\mathbf{GPD}$ we find the core (i.e. maximal sub-groupoid) of $\mathbf{Gpd}$ and $\mathcal{V}^{core}$. This allows us to define our cartesian and linear universes in any context $\Gamma$ as the functors:
\[
  \begin{split}
    \mathbb{U} :  \Gamma \to \mathbf{GPD}\\
    \mathbb{L} : \Gamma \to \mathbf{GPD}
  \end{split}
\]
which are constant at $\mathbf{Gpd}^{core}$ and $\mathcal{V}^{core}$, respectively. Any section $s : \Gamma \to \Gamma.\mathbb{U}$ will determine a functor $\hat s : \Gamma \to \mathbf{Gpd}$, which we  embed $\mathbf{Gpd} \to \mathbf{GPD}$ to get an interpretation of $El(s)$. Similarly, we get from each section $s : \Gamma \to \Gamma.\mathbb{L}$, a functor $El(s) : \Gamma \to \mathcal{V}$. 
It is easily seen that defining the linear universe via the core of $\mathcal{V}$ gives rise to the following interesting property, hinting at the possibility of a linear univalence axiom:
\begin{proposition}For a linear universe defined as above via $\mathcal{V}^{core}$, and two sections $s, t : \Gamma \to \Gamma.\mathbb{L}$, an isomorphism $\alpha : El(t) \cong El(s)$ gives rise to a section $p : \Gamma \to \Gamma.Id_{\mathbb{L}}\{s\}\{t\}$.
\end{proposition}
\subsection{Univalence in linear dependent types}\label{highermodel}
A key feature of the groupoid model is that it provides a model of dependent type theory where there might be nontrivial terms of the identity type. A natural question to ask is whether this higher dimensional feature of type theory can be extended to the linear dependent setting.
In particular, we might wish for a linear analogue to the univalence axiom to hold:
\[
  \inference{
    \Gamma \vdash A : \mathbb{L}\\
    \Gamma \vdash B : \mathbb{L}\\
    \Gamma; \cdot \vdash f : El(A) \multimap El(B)\\
    \Gamma; \cdot \vdash g : El(B) \multimap El(A)\\
    \Gamma; \cdot \vdash h : El(B) \multimap El(A)\\
    \Gamma \vdash p : \sigma(g \circ f)  =_{(El(A) \multimap El(A))_M} (id_A)_M\\
    \Gamma \vdash q : \sigma(f \circ h) =_{(El(B) \multimap El(B))_M} (id_B)_M\\
    }{\Gamma \vdash ua(f) : A =_{\mathbb{L}} B}[L-ua-I]
  \]
  To define the corresponding computation rule, we will the make use of a linear version of \texttt{transport}, which is easily definable through identity elimination. Given $\Gamma, x : C \vdash D \linear$, and an identity $p : a =_C b$, we get a function:
\[
  p^* : D[a/x] \multimap D[b/x]
  \]
  which we call the linear \textbf{transport along p}.
  This function has an inverse and thus yields for any $q : A =_L B$ a linear equivalence $El(A) \cong El(B)$.
  
    The computation rules for the univalence axiom asserts that the process creating equivalences from identities forms an inverse to univalence:
  \[
    \inference{
      \Gamma \vdash ua(f)^* : El(A) \multimap El(B)}
    {\Gamma \vdash ua(f)^* \equiv f : El(A) \multimap El(B)}
    [L-ua-C$_1$]
  \]
    \[
    \inference{
      \Gamma \vdash ua(p^*) : A =_\mathbb{L} B}
    {\Gamma \vdash ua(p^*) \equiv p : A =_{\mathbb{L}} B}
    [L-ua-C$_2$]
  \]
  The semantic interpretation of the procedure of turning an identity to an equivalence becomes the following:
  \begin{lemma}In the diagram model, equivalences $El(A) \cong El(B)$ are in one-to-one correspondence with identities $p : A =_L B$.
    \begin{proof}A section $p : \Gamma \to \Gamma.{\mathcal{V}^{core}}^{\to}$ defines for every morphism $\alpha : \gamma \to \gamma'  \in \Gamma$ a naturality square:
      \[
      \begin{tikzcd}
        A(\gamma) \ar[r,"p(\gamma)"] \ar[d,"A(\alpha)"] & B(\gamma) \ar[d,"B(\alpha)"]\\
        A(\gamma') \ar[r, "p(\gamma')"] & B(\gamma')
      \end{tikzcd}
      \]
where $A(\gamma)$ and $B(\gamma)$ are isomorphisms. Conversely, every natural isomorphism $El(A) \cong El(B)$ clearly defines such a section.
    \end{proof}
  \end{lemma}
  \begin{theorem}
    When $M$ factors through $\mathbf{Sets}$ as in Remark \ref{factorsets}, the linear univalence axiom holds in the diagram model. That is, given the following data:
    \begin{itemize}
    \item sections: $A, B : \Gamma \to \Gamma.{\mathbb{L}}$
    \item morphisms: $f : I \to [El(A), El(B)]$ and $g, h : I \to [El(B),El(A)]$ in $[\Gamma, \mathcal{V}]$,
    \item sections: $p : \Gamma \to \Gamma.(M(g f))^*(M(id_A))^*Id_{(M[El(A),El(A)])}$
    \item and $q : \Gamma \to \Gamma.(M(f g))^*(M(id_B))^*Id_{(M[El(B),El(B)])}$
    \end{itemize}
There is a natural isomorphism $El(A) \cong El(B)$.
    \begin{proof}
      The section $p$ selects for every $\gamma \in \Gamma$ an isomorphism between $\sigma(gf)$ and $\sigma(id_A)$ of $M[El(A),El(A)]$. But since $M$ factors through sets, $M[El(A), El(A)]$ is a discrete groupoid, so $\sigma(gf)$ and $\sigma(id_A)$ must be identical, and the same is true for $\sigma(fg)$ and $\sigma(id_B)$. Transporting back through the isomorphisms $[\Gamma,\mathbf{GPD}](1,M[El(A),El(A)]) \cong [\Gamma, \mathcal{V}](I, [El(A),El(A)]) \cong [\Gamma,\mathcal{V}](El(A),El(A))$ we find $gf = 1_{El(A)}$, and similarly $fg = 1_{El(B)}$.
\end{proof}
  \end{theorem}
  By the lemma above, equivalences $El(A) \cong El(B)$ are in one-to-one correspondence with identities $p : A =_L B$, demonstrating that the linear univalence axiom holds in the diagrams model as long as $M$ factors through sets.
  \subsubsection{Examples}
To summarize, these are the conditions imposed on $\mathcal{V}$ in order for the diagram model to support all of the type formers of our theory (including a universe of linear types satisfying the univalence axiom):
\begin{itemize}
\item A bicomplete symmetric monoidal closed category $\mathcal{V}$
\item An adjunction $L \dashv M$ between $\mathcal{V}$ and $\mathbf{Sets}$, such that $L\{*\}$ is isomorphic to the unit of the monoidal structure of $\mathcal{V}$.
\end{itemize}
Some concrete choices for $\mathcal{V}$ that fulfill these conditions are:
\begin{itemize}
\item The category $\mathbf{AbGroups}$ of abelian groups with the monoidal structure given by the tensor product of abelian groups. Here $L \dashv M$ arises from the free functor on abelian groups.
\item More generally, for any commutative ring $R$, the category $R$-$\mathbf{Mod}$ of modules over $R$ with the free functor/forgetful functor adjunction
\item The category $\mathbf{CGTop}_*$, of pointed compactly generated topological spaces, with the smash product as monoidal structure. The functor $M$ is here the forgetful functor which both forgets the base point and the topology, which has a left adjoint given by the discrete topology, and then taking the coproduct with the point to create a pointed space. The unit of $\mathbf{CGTop}_*$ is the two point discrete set $S^0$, which is precisely the image of the point in the adjunction above.
\end{itemize}
  \subsection{Discussion}
  Although the diagram model supports the univalence axiom, we are forced to truncate any higher dimensional structure by factoring $M$ through $\mathbf{Set}$, just as we can support it in the groupoid model for a universe only containing discrete groupoids. From the perspective of homotopy type theory, we may think of the set-indexed families model as a 0-dimensional model of linear dependent type theory and the diagram model as a 1-dimensional one. We conclude the paper by sketching what a 2-dimensional model might look:

  As outlined in \cite{smcat}, there is a symmetric monoidal structure on $\mathbf{SMCat}$, the category of small symmetric monoidal categories, symmetric monoidal functors and monoidal natural transformations. \footnote{Technically, the structure on $\mathbf{SMCat}$ is not quite symmetric monoidal, as the associators, unitors and symmetry functors are only invertible up to higher homotopy. However, if one applies these homotopies whenever necessary, one does get a model of linear dependent type theory.}
\begin{definition}
 Let the \textbf{2-dimensional model of LDTT} be given by the diagrams model where $\mathcal{V}$ is the 2-category of small symmetric monoidal categories, symmetric monoidal functors and monoidal natural transformations:
  \[
    \begin{tikzcd}
    \diag(\mathbf{SMCat}) \ar[rd, "cod"] & \diag(\mathbf{Gpd}) \ar[r, "\pi"] \ar[d, "cod"] & \mathbf{Gpd}^\to \ar[ld, "dom"]\\
    & \mathbf{Gpd}
    \end{tikzcd}
  \]
\end{definition}
For two symmetric monoidal categories $\mathcal{A}$ and $\mathcal{B}$, the (monoidal) functor category $[\mathcal{A}, \mathcal{B}]$ between them carries a natural monoidal structure \cite{smcat}, and gives $\mathbf{SMCat}$ a monoidal closed structure. Since $\mathbf{SMCat}$ is complete, with limits inherited from $\mathbf{Cat}$ equipped with a pointwise monoidal structure, we have support for $\sqcap$ and $\&$, and theorem \ref{diaghoms} gives us that this model supports $\multimap$ type formers. \footnote{Note, however, that we do not have all coproducts in $\mathbf{SMCat}$. Therefore, we cannot support $\oplus$ or $\sqsubset$. An alternative to be explored is the category $\mathbf{Mult}$, of multicategories, which is a symmetric monoidal closed, complete and co-complete \cite{elmendorf2009permutative}}
There is a natural candidate for the adjunction $L \dashv M$, based on the composite:
  \[
  \begin{tikzcd}
\mathbf{SMCat} \ar[r, "U"{name=D, below}, bend right]  & \ar[l, "F"{name=C, above}, bend right] \ar[from=C, to=D, symbol=\dashv] \mathbf{Cat} \ar[r, "\text{core}"{name=B, below}, bend right]   &  \mathbf{Gpd} \ar[l, "U"{name=A, above}, bend right] \ar[from=A, to=B, symbol=\dashv]
    \end{tikzcd}
  \]
  through which one should be able to construct a univalent universe containing nontrivial 1-dimensional linear types.

  Eventually, one would like to go all the way up and construct a $\infty$-dimensional formulation of linear dependent type theory. It has been speculated that models of a higher dimensional linear dependent type theory can be expressed through stable homotopy type theory \cite{schreiber2014quantization}, although it is unclear what the syntax for such a theory looks like.
\newpage
\bibliography{LundfallCSL}
\end{document}